\documentclass[10pt]{article}

\usepackage{amsmath}%
\usepackage{amsfonts}%
\usepackage{amssymb}%
\usepackage{multicol}%
\usepackage{graphicx}%
\usepackage{caption}
\usepackage{float}
\usepackage[titletoc]{appendix}

\newcommand{\hyp}[5]{\,\mbox{}_{#1}F_{#2}\!\left(\genfrac{}{}{0pt}{}{#3}{#4};#5\right)}

\numberwithin{equation}{section}

\usepackage{geometry}
\newtheorem{theorem}{Theorem}

\newtheorem{Remark}[theorem]{Remark}

\begin{document}

\title{Combinatorial identities and hypergeometric series}
\author{Enno Diekema \footnote{email adress: e.diekema@gmail.com}}
\maketitle

\begin{abstract}
\noindent
This paper describes a method to find a connection between combinatorial identities and hypergeometric series with a number of examples. Combinatorial identities can often be written as hypergeometric series with unit argument. In a number of cases these hypergeometric series are balanced and can be reduced to a simpler form. In this paper some combinatorial identities are proved using this method assuming that the results in the tables of Prudnikov et al. \cite{4} are proven without using hypergeometric functions.
\end{abstract}

\section{Introduction}
\setlength{\parindent}{0cm}
In a recent paper Qi et al. \cite{1} proved a number of combinatorial identities using the well-known series expansion of the function $f(x)=\arcsin(x)$. They claim that their method is better then other methods for example in the book \cite{2}. In addition, Qi also asks to prove another identity. 

An identity is called combinatorial if it contains at least one summation of one or more binomials. There is a vast literature on the subject and the most known books are \cite{13} and \cite{14}.

\

There are a number of methods of proving these identities. Following Gould \cite{13} we call the methods of comparison of coefficients in series expansions of a known function, use of differentiation or integration, finite difference methods, solution of recurrence formulas, mathematical induction, enumeration of lattice points, theory of combinations and permutations, and general series transformations. 

An application of the latter method is the use of hypergeometric functions. An example of this application is described in \cite[Ch.4]{8}. However, the discussion of this method is mainly limited to Gaussian hypergeometric functions. This paper  provides examples of the use of generalized hypergeometric functions. Bai-Ni Guo et.al. \cite{7} use this method to prove an identity of the Lah numbers. 

\

The advantage of this method compared to the combinatorial method, is that there is no need to search for some known function that can be expanded into a power series. The search for the well-known function is often a matter of luck or chance.

The disadvantage of this method is the fact that often some transformations of the hypergeometric function is needed. A great number of them are in \cite{4}. 

An additional advantage occurs if the hypergeometric function cannot be converted to a known form. If the summation is known by any method, then this gives a known form for the hypergeometric function. Then the table of the known hypergeometric functions, such as in [4] for example, can be significantly expanded. 

\

In this paper we use the following identities.

For the Pochhammer symbol
\[
(a)_n=\dfrac{\Gamma(a+n)}{\Gamma(a)}
\]
For the binomials
\begin{equation}
\binom{n}{k}=\dfrac{\Gamma(n+1)}{\Gamma)k+1)\Gamma(n+1-k)}=(-1)^k\dfrac{(-n)_k}{(1)_k}
\label{1.0}
\end{equation}
\begin{equation}
\binom{2k}{k}=\dfrac{\Gamma(2k+1)}{\Gamma(k+1)^2}=
\dfrac{2^{2k}}{\sqrt{\pi}}\dfrac{\Gamma\left(k+\dfrac{1}{2}\right)}{\Gamma(k+1)}
\qquad\quad
\label{1.1}
\end{equation}
The duplication formula for the Gamma function
\begin{equation}
\Gamma(2x)=\dfrac{1}{\sqrt{\pi}}\Gamma(x)\Gamma\left(x+\dfrac{1}{2}\right)2^{2x-1}
\label{1.3}
\end{equation}
The reflection formula for the Gamma function
\begin{equation}
\Gamma(x)\Gamma(1-x)=\dfrac{\pi}{\sin(x\, \pi)}
\label{1.3a}
\end{equation}
$\psi(x)$ is the digamma function and $\psi'(x)$ is the trigamma function. They are defined as
\[
\psi(z)=\dfrac{\Gamma'(x)}{\Gamma(x)} \qquad \text{and} \qquad
\psi'(x)=\dfrac{d}{dx}\dfrac{\Gamma'(x)}{\Gamma(x)}
\]
The trigamma function
\begin{equation}
\psi'(x)=\sum_{k=0}^\infty\dfrac{1}{(x+k)^2}
\label{1.4}
\end{equation}

It turns out that the used hypergeometric functions in many cases are terminated and zero or one-balanced. For the terminated one-balanced $_3F_2$ hypergeometrical  function with unit argument we use Saalsch\"utz's theorem
\begin{equation}
\hyp32{-n,a,b}{c,1+a+b-c-n}{1}=\dfrac{(c-a)_n(c-b)_n}{(c)_n(c-a-b)_n}
\label{1.7}
\end{equation}

\

\section{Overview of the combinatorial identities}
In this paper we treat the following identities. $S_0(n)$ to $S_3(n)$ are from \cite{1}. $S_4(n)$ was asked by Qi. $S_5(n)$ is from \cite{5}. $S_6(n)$ comes from exercise 2.45 from \cite{8}. $S_6(n)$ and $S_7(n)$ are identities from the internet. They were both proven with combinatorial methods \cite{9},\cite{12}. $S_8(n)$ and $S_9$ are from Gould \cite[1.79]{13} and \cite[2,23]{13}. They are added as an example of the method.

\begin{align}
&S_0(n)=\sum_{k=0}^n\dfrac{1}{(1+2k)}\binom{2k}{k}\binom{2(n-k)}{n-k}=
2^{4n}\dfrac{\Gamma(n+1)^2}{\Gamma(2n+2)}
\label{2.1}	\\
&S_1(n)=\sum_{k=0}^n\dfrac{1}{2^{4k}(n-k+1)^2}
\dfrac{\binom{2k}{k}}{\binom{2(n-k+1)}{n-k+1}}=
\dfrac{1}{2^{4n+3}}\dfrac{\Gamma(2n+2)}{\Gamma(n+1)\Gamma(n+2)}
\left[\dfrac{\pi^2}{2}-\psi'\left(\dfrac{3}{2}+n\right)\right]
\label{2.2}	\\
&S_2(n)=\sum_{k=0}^n\dfrac{1}{2^{4k}(2k+1)(n-k+1)}
\dfrac{\binom{2k}{k}}{\binom{2(n-k+1)}{n-k+1}}=
\dfrac{1}{2^{4n+3}}\dfrac{\Gamma(2n+2)}{\Gamma(n+1)\Gamma(n+2)}
\left[\dfrac{\pi^2}{2}-\psi'\left(\dfrac{3}{2}+n\right)\right]
\label{2.3} \\
&S_3(n)=\sum_{k=0}^n\dfrac{1}{2^{4k}(2k+1)(n-k+1)^2}
\dfrac{\binom{2k}{k}}{\binom{2(n-k+1)}{n-k+1}}=
\dfrac{3}{(2n+3)}\dfrac{1}{2^{4n+3}}\dfrac{\Gamma(2n+2)}{\Gamma(n+1)\Gamma(n+2)}
\left[\dfrac{\pi^2}{2}-\psi'\left(\dfrac{3}{2}+n\right)\right]
\label{2.4}	 \\
&S_4(n)=\sum_{k=0}^n\binom{2n+1}{2k}\binom{k}{m}=2^{2(n-m)}\dfrac{2n+1}{2(n-m)+1}
\binom{2n-m}{m} \qquad  m>1
\label{2.5} \\
&S_5=\sum_{k=1}^\infty\dfrac{(-1)^k}{k}\dfrac{1}{2^{2k}}\binom{2k}{k}=
2\ln\big[2\big(\sqrt{2}-1\big)\big]
\label{2.6} \\
&S_6(n)=L_{n+1,k+1}=\sum_{j=0}^n L_{j,k}\big(n+k+1\big)^{\overline{n-j}}
\label{2.7} \\
&\qquad\quad\ \  L_{i,j} \ \ \text{are the Lah numbers \cite{8}.} \nonumber \\
&\qquad\quad\ \ \big(n+k+1\big)^{\overline{n-j}} \ \text{ is the rising factorial of } (n+k+1) \nonumber \\
&S_7(n)=\sum_{k=1}^n\left[\dfrac{k}{n}\binom{2n}{n-k}\right]^2=C_{2n-1}
\label{2.8} \\
&\qquad\quad\ \  C_n \ \ \text{are the Catalan numbers \cite{11}.} \nonumber \\
&S_8(n)=\sum_{k=0}^n\binom{n+k}{k}\dfrac{1}{2^k}=2^n 
\label{2.9} \\
&S_9=\sum_{k=0}^\infty\dfrac{(-1)^k}{(2k+1)^2\binom{2k}{k}}=
\dfrac{\pi^2}{6}-3\ln^2\left(\dfrac{\sqrt{5}-1}{2}\right)
\label{2.10}
\end{align}

\

\section{Proofs of the identities}

\textbf{\fontsize{10.5}{12.5}\selectfont Proof of identity 2.1.}

Using \eqref{1.1} we get
\begin{align*}
S_0(n)
&=\sum_{k=0}^n\dfrac{1}{1+2k}\binom{2k}{k}\binom{2(n-k)}{n-k} \\
&=\sum_{k=0}^n\dfrac{\Gamma(1+2k)}{\Gamma(2+2k)}
\dfrac{\Gamma\left(\dfrac{1}{2}+k\right)2^{2k}}{\Gamma(1+k)\sqrt{\pi}}
\dfrac{\Gamma\left(\dfrac{1}{2}+n-k\right)2^{2(n-k)}}{\Gamma(1+n-k)\sqrt{\pi}}
\end{align*}
Using the duplication formula for the Gamma function in the first factor gives
\[
S_0(n)=\dfrac{2^{2n-1}}{\pi}
\sum_{k=0}^n
\dfrac{\Gamma\left(\dfrac{1}{2}+k\right)}{\Gamma\left(\dfrac{3}{2}+k\right)}
\dfrac{\Gamma\left(\dfrac{1}{2}+k\right)}{\Gamma(1+k)}
\dfrac{\Gamma\left(\dfrac{1}{2}+n-k\right)}{\Gamma(1+n-k)}
\]
After using $\Gamma(a+k)=\Gamma(a)(a)_k$ for the Gamma functions we get 
\[
S_0(n)=\dfrac{\left(\dfrac{1}{2}+n\right)}{\Gamma(1+n)}\dfrac{2^{2n}}{\sqrt{\pi}}
\sum_{k=0}^n\dfrac{(-n)_k\left(\dfrac{1}{2}\right)_k\left(\dfrac{1}{2}\right)_k}
{\left(\dfrac{3}{2}\right)_k\left(\dfrac{1}{2}-n\right)_k}\dfrac{1}{k!}
\]
Ths summation can be written as a $_3F_2$ hypergeometric function.
\[
S_0(n)=\dfrac{\left(\dfrac{1}{2}+n\right)}{\Gamma(1+n)}\dfrac{2^{2n}}{\sqrt{\pi}}
\hyp32{-n,\dfrac{1}{2},\dfrac{1}{2}}{\dfrac{3}{2},\dfrac{1}{2}-n}{1}
\]
The hypergeometric function is terminated and Saalsch\"utzian. Using \eqref{1.7} results in
\[
S_0(n)=\dfrac{\left(\dfrac{1}{2}+n\right)}{\Gamma(1+n)}\dfrac{2^{2n}}{\sqrt{\pi}}
\dfrac{(1)_n(1)_n}{\left(\dfrac{3}{2}\right)_n\left(\dfrac{1}{2}\right)_n}
\]
After some simplification we get \eqref{2.1} which proves the theorem. $\square$

\

\textbf{\fontsize{10.5}{12.5}\selectfont Proof of identity 2.2.}

Using \eqref{2.2} we get
\begin{align*}
S_1(n)&=\sum_{k=0}^n\dfrac{1}{2^{4k}(n-k+1)^2}
\dfrac{\binom{2k}{k}}{\binom{2(n-k+1)}{n-k+1}} \\
&=\dfrac{1}{2^{2n+2}}\sum_{k=0}^n\dfrac{\Gamma(n+1-k)^2}{\Gamma(n+2-k)}
\dfrac{\Gamma\left(\dfrac{1}{2}+k\right)}
{\Gamma(1+k)\Gamma\left(\dfrac{3}{2}+n-k\right)}
\end{align*}
Inversion of the summation with $m=n-k$ gives
\[
S_1(n)=\dfrac{1}{2^{2n+2}}\sum_{m=0}^n\dfrac{\Gamma(m+1)^2}{\Gamma(m+2)}
\dfrac{\Gamma\left(\dfrac{1}{2}+n-m\right)}
{\Gamma(n+1-m)\Gamma\left(\dfrac{3}{2}+m\right)}
\]
The Gamma functions can be written as Pochhammer symbols
\[
S_1(n)=\dfrac{\Gamma\left(\dfrac{1}{2}+n\right)}{\Gamma(1+n)\sqrt{\pi}}\dfrac{1}{2^{2n+1}}
\sum_{k=0}^n\dfrac{(-n)_m(1)_m(1)_m(1)_m}{(2)_m\left(\dfrac{3}{2}\right)_m\left(\dfrac{1}{2}-n\right)_m}\dfrac{1}{m!}
\]
The summation can be written as a hypergeometric function
\begin{align}
S_1(n)&=\dfrac{\Gamma\left(\dfrac{1}{2}+n\right)}{\Gamma(1+n)\sqrt{\pi}}\dfrac{1}{2^{2n+1}}
\hyp43{-n,1,1,1}{2,\dfrac{3}{2},\dfrac{3}{2},\dfrac{1}{2}-n}{1} \nonumber \\
&=\dfrac{1}{2^{4n+1}}\dfrac{\Gamma(2n+1)}{\Gamma(n+1)^2}
\hyp43{-n,1,1,1}{\dfrac{3}{2},\dfrac{1}{2}-n,2}{1}
\label{3.1}
\end{align}
The hypergeometric function is Saalsch\"utzian. There are many transformations of this hypergeometric function. We use \cite[(1.6)]{3}
\[
\hyp43{-n,a_1,a_2,a_3}{b_1,b_2,1-s-n}{1}=
\dfrac{(a_1+s)_n(a_2+s)_n(a_3)_n}{(b_1)_n(b_2)_n(s)_n}
\hyp43{b_1-a_3,b_2-a_3,s,-n}{a_1+s,a_2+s,1-a_3-n}{1}
\]
with $s=b_1+b_2-a_1-a_2-a_3$. 

Application to the hypergeometric function with $a_1=a_2=a_3=1,\ b_1=2,\ b_2=\dfrac{3}{2}$ gives
\[
\hyp43{-n,1,1,1}{2,\dfrac{3}{2},\dfrac{1}{2}-n}{1}=
\dfrac{\left(\dfrac{3}{2}\right)_n(1)_n}{\left(\dfrac{1}{2}\right)_n(2)_n}
\hyp43{1,\dfrac{1}{2},\dfrac{1}{2},-n}{\dfrac{3}{2},\dfrac{3}{2},-n}{1}
\]
Simplification gives
\begin{align}
\hyp43{-n,1,1,1}{2,\dfrac{3}{2},\dfrac{1}{2}-n}{1}
&=\dfrac{(2n+1)}{(n+1)}\hyp43{1,\dfrac{1}{2},\dfrac{1}{2},-n}{\dfrac{3}{2},\dfrac{3}{2},-n}{1}= 
\dfrac{(2n+1)}{(n+1)}\sum_{k=0}^n\dfrac{(1)_k\left(\dfrac{1}{2}\right)_k\left(\dfrac{1}{2}\right)_k(-n)_k}{\left(\dfrac{3}{2}\right)_k\left(\dfrac{3}{2}\right)_k(-n)_k}
\dfrac{1}{k!}= \nonumber \\
&=\dfrac{(2n+1)}{4(n+1)}
\dfrac{\Gamma\left(\dfrac{3}{2}\right)^2}{\Gamma\left(\dfrac{1}{2}\right)^2}
\sum_{k=0}^n
\dfrac{\Gamma\left(k+\dfrac{1}{2}\right)^2}{\Gamma\left(k+\dfrac{3}{2}\right)^2}=
\dfrac{(2n+1)}{4(n+1)}
\sum_{k=0}^n\dfrac{1}{\left(k+\dfrac{1}{2}\right)^2}= \nonumber \\
&=\dfrac{(2n+1)}{4(n+1)}
\sum_{k=0}^\infty\dfrac{1}{\left(k+\dfrac{1}{2}\right)^2}-
\dfrac{(2n+1)}{4(n+1)}
\sum_{k=n+1}^\infty\dfrac{1}{\left(k+\dfrac{1}{2}\right)^2}= \nonumber \\
&=
\dfrac{(2n+1)}{4(n+1)}\left(\dfrac{\pi^2}{2}-\psi'\left(n+\dfrac{3}{2}\right)\right)
\label{3.2}
\end{align}
Substitution in \eqref{3.1} results in
\[
S_1=\dfrac{1}{2^{4n+1}}\dfrac{\Gamma(2n+1)}{\Gamma(n+1)^2}
\dfrac{(2n+1)}{4(n+1)}\left(\dfrac{\pi^2}{2}-\psi'\left(n+\dfrac{3}{2}\right)\right)
\]
Simplification proves the theorem. $\square$

\

\textbf{\fontsize{10.5}{12.5}\selectfont Proof of identity 2.3.}

Using \eqref{2.3} we get
\begin{align*}
S_2(n)&=\sum_{k=0}^n\dfrac{1}{2^{4k}(2k+1)(n-k+1)}
\dfrac{\binom{2k}{k}}{\binom{2(n-k+1)}{n-k+1}} \\
&=\dfrac{1}{2^{2n+3}}\sum_{k=0}^n\dfrac{\Gamma(n+1-k)\Gamma\left(\dfrac{1}{2}+k\right)^2}
{\Gamma\left(\dfrac{3}{2}+k\right)\Gamma(1+k)\Gamma\left(\dfrac{3}{2}+n-k\right)}
\end{align*}
Inversion of the summation with $m=n-k$ gives
\begin{align}
S_2(n)&=\dfrac{1}{2^{2n+3}}\sum_{m=0}^n\dfrac{\Gamma(m+1)\Gamma\left(\dfrac{1}{2}+n-m\right)^2}
{\Gamma\left(\dfrac{3}{2}+n-m\right)\Gamma(1+n-m)\Gamma\left(\dfrac{3}{2}+m\right)} \nonumber \\
&=\dfrac{1}{2^{2n+3}}\dfrac{\Gamma\left(\dfrac{1}{2}+n\right)^2}
{\Gamma\left(\dfrac{3}{2}+n\right)\Gamma\left(\dfrac{3}{2}\right)\Gamma(n+1)}
\sum_{m=0}^n\dfrac{(-n)_m(1)_m(1)_m\left(-\dfrac{1}{2}-n\right)_m}
{\left(\dfrac{1}{2}-n\right)_m\left(\dfrac{1}{2}-n\right)_m\left(\dfrac{3}{2}\right)_m}
\dfrac{1}{m!} \nonumber \\
&=\dfrac{1}{2^{2n+1}}\dfrac{\left(\dfrac{1}{2}\right)_n}
{(1+2n)(1)_n}\hyp43{-n,1,1,-\dfrac{1}{2}-n}{\dfrac{3}{2},\dfrac{1}{2}-n,\dfrac{1}{2}-n}{1}
\label{3.3}
\end{align}
It is an easy task to prove by mathematical induction the following property
\[
\hyp43{-n,1,1,-\dfrac{1}{2}-n}{\dfrac{3}{2},\dfrac{1}{2}-n,\dfrac{1}{2}-n}{1}=
(2n+1)\hyp43{-n,1,1,1}{\dfrac{3}{2},2,\dfrac{1}{2}-n}{1}
\]
Substitution of the result from \eqref{3.2} gives
\begin{equation}
\hyp43{-n,1,1,-\dfrac{1}{2}-n}{\dfrac{3}{2},\dfrac{1}{2}-n,\dfrac{1}{2}-n}{1}=
\dfrac{(2n+1)^2}{4(n+1)}\left(\dfrac{\pi^2}{2}-\psi'\left(n+\dfrac{3}{2}\right)\right)
\label{3.4}
\end{equation}
Substitution in \eqref{3.3} results in
\[
S_2(n)=\dfrac{1}{2^{2n+1}}\dfrac{\left(\dfrac{1}{2}\right)_n}
{(1+2n)(1)_n}
\dfrac{(2n+1)^2}{4(n+1)}\left(\dfrac{\pi^2}{2}-\psi'\left(n+\dfrac{3}{2}\right)\right)
\]
After some manipulation with the Pochhammer symbols this completes the proof. $\square$

\

\textbf{\fontsize{10.5}{12.5}\selectfont Proof of identity 2.4.}

Using \eqref{2.4} we get
\begin{align*}
S_3(n)&=\sum_{k=0}^n\dfrac{1}{2^{4k}(2k+1)(n-k+1)^2}
\dfrac{\binom{2k}{k}}{\binom{2(n-k+1)}{n-k+1}} \\
&=\dfrac{1}{2^{2n+3}}\sum_{k=0}^n\dfrac{\Gamma(n+1-k)^2\Gamma\left(\dfrac{1}{2}+k\right)^2}
{\Gamma\left(\dfrac{3}{2}+k\right)\Gamma(n+2-k)\Gamma(1+k)\Gamma\left(\dfrac{3}{2}+n-k\right)}
\end{align*}
Inversion of the summation with $m=n-k$ gives
\begin{align}
S_3(n)&=\dfrac{1}{2^{2n+3}}\sum_{m=0}^n\dfrac{\Gamma(m+1)^2\Gamma\left(\dfrac{1}{2}+n-m\right)^2}
{\Gamma\left(\dfrac{3}{2}+n-m\right)\Gamma(m+2)\Gamma(1+n-m)\Gamma\left(\dfrac{3}{2}+m\right)} \nonumber \\
&=\dfrac{1}{2^{2n+3}}\dfrac{\Gamma\left(\dfrac{1}{2}+n\right)^2}
{\Gamma\left(\dfrac{3}{2}+n\right)\Gamma\left(\dfrac{3}{2}\right)\Gamma(n+1)}
\sum_{m=0}^n\dfrac{(-n)_m(1)_m(1)_m(1)_m\left(-\dfrac{1}{2}-n\right)_m}
{(2)_m\left(\dfrac{1}{2}-m\right)\left(\dfrac{1}{2}-m\right)\left(\dfrac{3}{2}\right)}
\dfrac{1}{m!} \nonumber \\
&=\dfrac{1}{2^{2n+3}}\dfrac{\left(\dfrac{1}{2}\right)_n}
{(1+2n)(1)_n}\hyp54{-n,1,1,1,-\dfrac{1}{2}-n}{\dfrac{3}{2},\dfrac{1}{2}-n,2,\dfrac{1}{2}-n}{1}
\label{3.5}
\end{align}
The hypergeometric function is two-balanced. Application of \cite[7.2.3.(20)]{4} with $\rho=1 \ \ \text{and}\ \ \sigma=-\dfrac{1}{2}-n$ \ gives
\[
\hyp54{-n,1,1,1,-\dfrac{1}{2}-n}{\dfrac{3}{2},\dfrac{1}{2}-n,2,\dfrac{1}{2}-n}{1}=
\dfrac{2n+1}{2n+3}\hyp43{-n,1,1,1}{\dfrac{3}{2},\dfrac{1}{2}-n,2}{1}+
\dfrac{2}{3+2n}\hyp43{-n,1,1,-\dfrac{1}{2}-n}{\dfrac{3}{2},\dfrac{1}{2}-n,\dfrac{1}{2}-n}{1}
\]
The first hypergeometric function in the right-hand side has already computed \eqref{3.2}. The second hypergeometric function in the right-hand side has also computed \eqref{3.4}. Then we get
\[
\hyp54{-n,1,1,1,-\dfrac{1}{2}-n}{\dfrac{3}{2},\dfrac{1}{2}-n,2,\dfrac{1}{2}-n}{1}=
\dfrac{3}{4}\dfrac{(2n+1)^2}{(3+2n)(n+1)}\left(\dfrac{\pi^2}{2}-\psi'\left(n+\dfrac{3}{2}\right)\right)
\]
Substitution in \eqref{3.5} proves the theorem. $\square$

\

\textbf{\fontsize{10.5}{12.5}\selectfont Proof of identity 2.5.}

Using \eqref{2.5} we get
\[
S_4(n)=\sum_{k=0}^n\binom{2n+1}{2k}\binom{k}{m}=\dfrac{\Gamma(2n+2)}{\Gamma(m+1)}
\sum_{k=0}^n\dfrac{1}{\Gamma(1+2k)\Gamma(2n+2-2k)}
\dfrac{\Gamma(1+k)}{ \Gamma(1-m+k)}
\]
Application of \eqref{1.3} gives
\[
S_4(n)=\dfrac{\pi}{2^{2n+1}}\dfrac{\Gamma(2n+2)}{\Gamma(m+1)}
\sum_{k=0}^n\dfrac{1}{\Gamma\left(\dfrac{1}{2}+k\right)\Gamma(n+1-k)\Gamma\left(n+\dfrac{3}{2}-k\right)}\dfrac{1}{\Gamma(1-m+k)}
\]
After using $\Gamma(a-k)=(-1)^n\dfrac{\Gamma(a)}{(1-a)_k}$ for the Gamma functions we get
\begin{align*}
S_4(n)&=\dfrac{\pi}{2^{2n+1}}\dfrac{\Gamma(2n+2)}{\Gamma(m+1)}
\dfrac{1}{\left(\dfrac{1}{2}\right)\Gamma(n+1)\Gamma\left(n+\dfrac{3}{2}\right)\Gamma(1-m)}
\sum_{k=0}^n\dfrac{(-n)_k\left(-\dfrac{1}{2}-n\right)_k(1)_k}{\left(\dfrac{1}{2}\right)_k(1-m)_k}\dfrac{1}{k!} \\
&=\dfrac{1}{\Gamma(1+m)\Gamma(1-m)}\hyp32{-n,\dfrac{1}{2},1}{1-m,\dfrac{1}{2}}{1}
\end{align*}
Because one of the lower parameters is negative ($1-m<0$) we use \cite[7.2.3.(6)]{4}
which we write in a modified form.
\begin{multline*}
\dfrac{1}{\Gamma(-M)} \ _{p+1}F_p
\left(\begin{array}{l}
	a_0,\dots,a_p \\
	-M,b_2,\dots,b_p
\end{array};x\right)= \\
=\dfrac{x^{M+1}(a_0)_{M+1}\dots(a_p)_{M+1}}{\Gamma(M+2)(b_2)_{M+1}\dots(b_p)_{M+1}}
\ _{p+1}F_p
\left(\begin{array}{l}
	a_0+M+1,\dots,a_p+M+1 \\
	M+2,b_2+M+1,\dots,b_p+M+1
\end{array};x\right)
\end{multline*}
with $M$ a positive integer. For the proof of this formula see \cite[Lemma 2]{6}. Application gives
\[
S_4(n)=\dfrac{1}{\Gamma(m+1)}\dfrac{(-n)_m\left(-\dfrac{1}{2}-n\right)_m}
{\left(\dfrac{1}{2}\right)_m}\hyp21{m-n,m-\dfrac{1}{2}-n}{m+\dfrac{1}{2}}{1}
\]
The hypergeometric function can be summed.
\[
S_4(n)=\dfrac{1}{\Gamma(m+1)}\dfrac{(-n)_m\left(-\dfrac{1}{2}-n\right)_m}
{\left(\dfrac{1}{2}\right)_m}
\dfrac{\Gamma\left(m+\dfrac{1}{2}\right)\Gamma(2n-m+1)}{\Gamma\left(n+\dfrac{1}{2}\right)\Gamma(n+1)}
\]
Using $(-n)_m=(-1)^m\frac{\Gamma(n+1)}{\Gamma(n-m)}$ and elaboration of the binomials results in
\begin{align*}
S_4(n)&=(-1)^m\dfrac{\Gamma(n+1)\Gamma\left(m-\dfrac{1}{2}-n\right)}
{\Gamma(n-m+1)\Gamma\left(-\dfrac{1}{2}-n\right)\Gamma(2n+1)}
\dfrac{\Gamma(2n-m+1)2^{2n}}{\Gamma(m+1)} \\
&=(-1)^m\dfrac{\Gamma\left(n+\dfrac{3}{2}\right)\Gamma(n+1)\Gamma\left(n+\dfrac{3}{2}-m\right)\Gamma\left(m-\dfrac{1}{2}-n\right)}
{\Gamma\left(n+\dfrac{3}{2}\right)\Gamma\left(-\dfrac{1}{2}-n\right)\Gamma(n-m+1)\Gamma\left(n-m+\dfrac{3}{2}\right)\Gamma(2n+1)}
\dfrac{\Gamma(2n-m+1)2^{2n}}{\Gamma(m+1)}
\end{align*}
Using the reflection formula \eqref{1.3a} and the duplication formula \eqref{1.3} gives at last
\[
S_4(n)=\dfrac{(2n+1)\Gamma(2n-m+1)2^{2n-2m}}{\Gamma(2n-2m+2)\Gamma(1+m)}
\]
The right-hand side of \eqref{2.5} gives
\[
2^{2(n-m)}\dfrac{2n+1}{2(n-m)+1}\binom{2n-m}{m}=
2^{2(n-m)}\dfrac{(2n+1)}{(2n-2m+1)}\dfrac{\Gamma(2n-m+1)}{\Gamma(m+1)\Gamma(2n-2m+1)}
\]
Simplifying proves the property. $\square$

\

\textbf{\fontsize{10.5}{12.5}\selectfont Proof of identity 2.6.}

Using \eqref{2.6} and \eqref{1.1} we get
\[
S_5=\sum_{k=1}^\infty\dfrac{(-1)^k}{k}\dfrac{1}{2^{2k}}\binom{2k}{k}=
\sum_{k=1}^\infty(-1)^k\dfrac{2^{2k}\Gamma\left(k+\dfrac{1}{2}\right)}
{\sqrt{\pi}\Gamma(k+1)}\dfrac{1}{2^{2k}k}
\]
Using $k=m+1$ gives
\begin{align*}
S_5&=-\sum_{m=0}^\infty(-1)^m\dfrac{\Gamma\left(m+\dfrac{3}{2}\right)}
{\sqrt{\pi}\Gamma(m+2)(m+1)}=
-\dfrac{\Gamma\left(\dfrac{3}{2}\right)}{\sqrt{\pi}}
\sum_{m=0}^\infty(-1)^m\dfrac{\left(\dfrac{3}{2}\right)_m\Gamma(m+1)}
{(m+1)\Gamma(m+1)\Gamma(m+2)}= \\
&=-\dfrac{1}{2}
\sum_{m=0}^\infty\dfrac{\left(\dfrac{3}{2}\right)_m(1)_m(1)_m}
{(2)_m(2)_m}\dfrac{1}{m!}(-1)^m=-\dfrac{1}{2}\hyp32{1,1,\dfrac{3}{2}}{2,2}{-1}
\end{align*}
\cite[7.4.1(365)]{4} gives
\[
\hyp32{1,1,\dfrac{3}{2}}{2,2}{z}=\dfrac{4}{z}\ln\left(\dfrac{2\big(1-\sqrt{1-z}\big)}{z}\right)
\]
Substitution of $z=-1$ proves the theorem. $\square$

\

\textbf{\fontsize{10.5}{12.5}\selectfont Proof of identity 2.7.}

We repeat \eqref{2.7}
\begin{equation}
S_6(n)=L_{n+1,k+1}=\sum_{j=0}^n L_{j,k}\big(n+k+1\big)^{\overline{n-j}}
\label{7.1}
\end{equation}
$L_{n,k}$ is the Lah number. It is defined as
\[
L_{n,k}=\binom{n-1}{k-1}\dfrac{n!}{k!}
\]
with $L_{0,0}=1$ and $L_{n,k}=0$ if $k>n$. So we can write for the left-hand side of \eqref{7.1}
\begin{equation}
L_{n+1,k+1}=\binom{n}{k}\dfrac{(n+1)!}{(k+1)!}=
\dfrac{\Gamma(n+1)\Gamma(n+2)}{\Gamma(k+1)\Gamma(k+2)\Gamma(n+1-k)}
\label{7.2}
\end{equation}
For the right-hand side of \eqref{7.1} we have $\big(n+k+1\big)^{\overline{n-j}}$ which is the rising factorial of $(n+k+1)$. We get
\[
\big(n+k+1\big)^{\overline{n-j}}=\dfrac{\Gamma(n+k+2)}{\Gamma(n+k+2-(n-j))}
=\dfrac{(k+2)_n}{(k+2)_j}
\]
For the right-hand side of \eqref{7.1} we get (because $0 \leq k\leq j\leq n$)
\[
R=\sum_{j=0}^n L_{j,k}\big(n+k+1\big)^{\overline{n-j}}=
\sum_{j=k}^n\dfrac{\Gamma(j)}{\Gamma(k)\Gamma(j-k+1)}\dfrac{j!}{k!}
\dfrac{(k+2)_n}{(k+2)_j}
\]
For the lower bound of the summation we substitute \  $j=m+k$
\begin{align*}
R&=\sum_{m=0}^{n-k}\dfrac{\Gamma(k+m)}{\Gamma(k)\Gamma(m+1)}\dfrac{(m+k)!}{k!}
\dfrac{(k+2)_n}{(k+2)_{k+m}} 
=\dfrac{(k+2)_n}{(k+2)_{k}}\sum_{m=0}^{n-k}\dfrac{(k)_m(k+1)_m}{(2k+2)_m}\dfrac{1}{m!}= \\
&=\dfrac{(k+2)_n}{(k+2)_{k}}\sum_{m=0}^\infty\dfrac{(k)_m(k+1)_m}{(2k+2)_m}\dfrac{1}{m!}-\dfrac{(k+2)_n}{(k+2)_{k}}\sum_{m=n-k+1}^\infty\dfrac{(k)_m(k+1)_m}{(2k+2)_m}\dfrac{1}{m!}
\end{align*}
The first summation is a hypergeometric function with unit argument and can be summed. For the second summation we substitute \ $m=j+n-k+1$.
\[
R=\dfrac{(k+2)_n}{(k+2)_{k}}\binom{2k+1}{k+1}-
\dfrac{(k+2)_n}{(k+2)_{k}}
\sum_{j=0}^\infty\dfrac{(k)_{j+n-k+1}(k+1)_{j+n-k+1}}{(2k+2)_{(j+n-k+1)}}
\dfrac{1}{(j+n-k+1)!}
\]
Simplifying gives
\[
R=\dfrac{\Gamma(k+n+2)}{\Gamma(k+1)\Gamma(k+2)}-
\dfrac{\Gamma(k+n+2)}{\Gamma(k)\Gamma(k+1)}
\dfrac{\Gamma(n+1)\Gamma(n+2)}{\Gamma(n+3+k)\Gamma(n+2-k)}
\sum_{j=0}^\infty\dfrac{(n+1)_j(n+2)_j}{(n+3+k)_j(n+2-k)_j}
\]
The summation can be written as a hypergeometric function
\[
R=\dfrac{\Gamma(k+n+2)}{\Gamma(k+1)\Gamma(k+2)}-
\dfrac{\Gamma(k+n+2)}{\Gamma(k)\Gamma(k+1)}
\dfrac{\Gamma(n+1)\Gamma(n+2)}{\Gamma(n+3+k)\Gamma(n+2-k)}
\hyp32{n+1,n+2,1}{n+k+3,n-k+2}{1}
\]
Transforming this hypergeometric function with \cite[16.4.11]{10} gives
\[
R=\dfrac{\Gamma(k+n+2)}{\Gamma(k+1)\Gamma(k+2)}-
\dfrac{\Gamma(n+1)\Gamma(n+2)}{\Gamma(k)\Gamma(k+1)\Gamma(n+2-k)}
\hyp32{1-k,-k,1}{n-k+2,2}{1}
\]
The \ $_3F_2$ \ hypergeometric function can be transformed into a \ $_2F_1$ \ hypergeometric function with \cite[7.2.3(17)]{4}
\[
_pF_q
\left(\begin{array}{l}
	(a_{p-1}),1 \\
	(b_{q-1}),2
\end{array};z\right)=\prod_{j=1}^{q-1}(b_j-1)\prod_{k=1}^{p-1}(a_k-1)^{-1}
\left[_{p-1}F_{q-1}
\left(\begin{array}{l}
	(a_{p-1})-1 \\
	(b_{q-1})-1
\end{array};z\right)-1 \right]
\]
Application gives
\[
R=\dfrac{\Gamma(k+n+2)}{\Gamma(k+1)\Gamma(k+2)}-
\dfrac{\Gamma(n+1)\Gamma(n+2)}{\Gamma(k+1)\Gamma(k+2)\Gamma(n-k+1)}
\left[\hyp21{-k,-1-k}{n-k+1}{1}-1\right]
\]
The hypergeometric function is known.
\[
\hyp21{-k,-1-k}{n-k+1}{1}=\dfrac{\Gamma(n-k+1)\Gamma(n+k+2)}{\Gamma(n+1)\Gamma(n+2)}
\]
Substitution gives after some simplification
\[
R=\dfrac{\Gamma(n+1)\Gamma(n+2)}{\Gamma(k+1)\Gamma(k+2)\Gamma(n-k+1)}
\]
Comparing with \eqref{7.2} proves the theorem. $\square$

\

\textbf{\fontsize{10.5}{12.5}\selectfont Proof of identity 2.8.}

We repeat \eqref{2.8}
\begin{equation}
S_7(n)=\sum_{k=1}^n\left[\dfrac{k}{n}\binom{2n}{n-k}\right]^2=C_{2n-1}
\label{8.1}
\end{equation}
The Catalan number is defined as
\[
C_n=\dfrac{\Gamma(2n+1)}{\Gamma(n+2)\Gamma(n+1)}
\]
which gives
\[
C_{2n-1}=\dfrac{\Gamma(4n-1)}{\Gamma(2n+1)\Gamma(2n)}
\]
The summation can be written as a hypergeometric function.
\begin{align}
\sum_{k=1}^n\left[\dfrac{k}{n}\binom{2n}{n-k}\right]^2
&=\sum_{k=1}^n\left[\dfrac{k}{n}\dfrac{\Gamma(2n+1)}{\Gamma(n+1-k)\Gamma(n+1+k)}
\right]^2 
=\dfrac{\Gamma(2n+1)^2}{n^2\Gamma(n+1)^4}\sum_{k=1}^n\left[\dfrac{(-n)_k}{(n+1)_k}
\dfrac{(1)_k}{(k-1)!}\right]^2 \nonumber \\
&=\dfrac{\Gamma(2n+1)^2}{n^2\Gamma(n+1)^4}\sum_{k=0}^{n-1}\left[\dfrac{(-n)_{k+1}}{(n+1)_{k+1}}\dfrac{(1)_{k+1}}{k!}\right]^2
=\dfrac{\Gamma(2n+1)^2}{\Gamma(n+1)^4(n+1)^2}\sum_{k=0}^{n-1}\left[\dfrac{(1-n)_k}{(n+2)_k}\dfrac{(2)_k}{k!}\right]^2 \nonumber \\
&=\dfrac{\Gamma(2n+1)^2}{\Gamma(n+2)^2\Gamma(n+1)^2}
\hyp43{2,2,1-n,1-n}{1,2+n,2+n}{1}
\label{8.1a}
\end{align}
The hypergeometric function is not directly known in closed form. Making use of a number of transformations we get a usable form. We start with \cite[7.2.3.(25)]{4} 
\[
\sigma\  _p F_q\left(
\begin{array}{l}
	(a_{p-1}),\rho,\sigma+1 \\
	\qquad (b_q)
\end{array};z \right)-
\rho\ _p F_q\left(
\begin{array}{l}
	(a_{p-1}),\rho+1,\sigma \\
	\qquad (b_q)
\end{array};z \right)=
(\sigma-\rho)\ _p F_q\left(
\begin{array}{l}
	(a_{p-1}),\rho,\sigma \\
	\qquad (b_q)
\end{array};z \right)
\]
Application with $\sigma=1$ and $\rho=2$ gives
\begin{equation}
\hyp43{2,2,1-n,1-n}{1,2+n,2+n}{1}=2\ \hyp32{1-n,1-n,3}{2+n,2+n}{1}-
\hyp43{1,1-n,1-n,2}{2+n,2+n,1}{1}
\label{8.2}
\end{equation}
For the second hypergeometric function on the right-hand side we use \cite[7.5.3.(6)]{4}
\begin{align*}
\hyp43{1,a,b,c}{3-a,3-b,3-c}{1}
&=\dfrac{1}{2(a-1)(b-1)(c-1)}
\Gamma\left[
\begin{array}{l}
	3-a,3-b,3-c,4-a-b-c \\
	3-a-b,3-a-c,3-b-c
\end{array}
\right]- \\
&-\dfrac{(2-a)(2-b)(2-c)}{2(a-1)(b-1)(c-1)}
\end{align*}
Application with $a=b=1-n$ and $c=2$ gives
\begin{equation}
\hyp43{1,1-n,1-n,2}{2+n,2+n,1}{1}=\hyp32{1-n,1-n,2}{2+n,2+n}{1}=\dfrac{(n+1)^2}{4n}
\label{8.3}
\end{equation}
For the first hypergeometric function on the right-hand side we use \cite[16.3.7]{10} with $z=1$
\begin{multline*}
\hyp32{a_1+1,a_2,a_3}{b_1,b_2}{1}a_1(b_1+b_2-a_1-a_2-a_2-1)+ \\
+\hyp32{a_1,a_2,a_3}{b_1,b_2}{1}\big((2a_1-b_1)(2a_1-b_2)+a_1-(a_1)^2-(a_1-a_2)(a_1-a_3) \big)- \\
-\hyp32{a_1-1,a_2,a_3}{b_1,b_2}{1}(a_1-b_1)(a_1-b_2)=0
\end{multline*}
With $a_1=2$, \ $a_2=a_3=1-n$ and $b_1=b_2=2+n$  we get
\[
\hyp32{3,1-n,1-n}{2+n,2+n}{1}2(4n-1)+\hyp32{1-n,1-n,2}{2+n,2+n}{1}(1-6n)-
\hyp32{1,1-n,1-n}{2+n,2+n}{1}=0
\]
Application of \eqref{8.3} gives
\begin{equation}
\hyp32{3,1-n,1-n}{2+n,2+n}{1}=\dfrac{n^2}{2(4n-1)}\hyp32{1,1-n,1-n}{2+n,2+n}{1}-
\dfrac{(n+1)^2}{8n}(1-6n)
\label{8.4}
\end{equation}
Substitution of \eqref{8.3} and \eqref{8.4} into \eqref{8.2} gives
\begin{equation}
\hyp43{1-n,1-n,2,2}{2+n,2+n,1}{1}=\dfrac{n^2}{(4n-1)}\hyp32{1-n,1-n,1}{2+n,2+n}{1}+
\dfrac{(n+1)^2}{2(4n-1)}
\label{8.5}
\end{equation}
The hypergeometric function on the right-hand side is not known. So we need still another transformation. Writing the hypergeometric function as a summation we get
\[
\hyp32{1-n,1-n,1}{2+n,2+n}{1}=\sum_{k=0}^{n-1}\dfrac{(1-n)_k(1-n)_k}{(2+n)_k(2+n)_k}
\]
The summation can be inverted with $m=(n-1)-k$. The result is
\[
\hyp32{1-n,1-n,1}{2+n,2+n}{1}=\dfrac{\Gamma(n)^2\Gamma(n+2)^2}{\Gamma(2n+1)^2}
\sum_{m=0}^{n-1}\dfrac{(-2n)_m(-2n)_m}{(1)_m}\dfrac{1}{m!}
\]
This summation cannot be written as a hypergeometric function, because the upper bound is $n-1$. So we split the summation
\[
\hyp32{1-n,1-n,1}{2+n,2+n}{1}=\dfrac{\Gamma(n)^2\Gamma(n+2)^2}{\Gamma(2n+1)^2}
\sum_{m=0}^\infty\dfrac{(-2n)_m(-2n)_m}{(1)_m}\dfrac{1}{m!}-
\sum_{m=n}^{2n}\dfrac{(-2n)_m(-2n)_m}{(1)_m}\dfrac{1}{m!}
\]
The first summation is a hypergeometric function. The second summation can be converted into a usable summation
\begin{multline*}
\hyp32{1-n,1-n,1}{2+n,2+n}{1}=
\dfrac{\Gamma(n)^2\Gamma(n+2)^2}{\Gamma(2n+1)^2}\hyp21{-2n,-2n}{1}{1}- \\
-\dfrac{\Gamma(n)^2\Gamma(n+2)^2}{\Gamma(2n+1)^2}\dfrac{(-2n)_n(-2n)_n}{(1)_n(1)_n}
\sum_{k=0}^n\dfrac{(-n)_k(-n)_k(1)_k}{(n+1)_k(n+1)_k}\dfrac{1}{k!}
\end{multline*}
The first hypergeometric function is well known. The summation can be written as a hypergeometric function. At last we get
\begin{equation}
\hyp32{1-n,1-n,1}{2+n,2+n}{1}=\dfrac{\Gamma(n)^2\Gamma(n+2)^2\Gamma(4n+1)}{\Gamma(2n+1)^4}-\dfrac{\Gamma(n)^2\Gamma(n+2)^2}{\Gamma(n+1)^4}
\hyp32{-n,-n,1}{1+n,1+n}{1}
\label{8.6}
\end{equation}
Now we can use \cite[7.4.4.(31)]{4}
\begin{multline*}
\hyp32{a,b,1}{-m-a,-m-b}{1}=\dfrac{1}{2}\sum_{k=0}^{m+1}\dfrac{(a)_k(b)_k}{(-a-m)_k(-b-m)_k}+ \\
+\dfrac{2^{2m-1}\sqrt{\pi}}{(1+m+2a)_m(1+m+2b)_m}
\Gamma\left[
\begin{array}{l}
	1-a,1-b,-a-b-m-\dfrac{1}{2}\\
	-a-b-m,\dfrac{1}{2}-a-m,\dfrac{1}{2}-b-m
\end{array}
\right]
\end{multline*}
With $a=b=-n$ and $m=-1$ we get after a lot of manipulations with the Gamma functions and the Poch- hammer symbols
\begin{equation}
\hyp32{-n,-n,1}{1+n,1+n}{1}=\dfrac{1}{2}+
\dfrac{1}{2}\dfrac{\Gamma(n+1)^4\Gamma(4n+1)}{\Gamma(2n+1)^4}
\label{8.7}
\end{equation}
This can also be written as
\[
\hyp32{-n,-n,1}{1+n,1+n}{1}=\dfrac{1}{2}+\dfrac{1}{2}\dfrac{\binom{4n}{2n}}{\binom{2n}{n}^2}
\]
Substitution \eqref{8.7} into \eqref{8.6} and the result into \eqref{8.5} gives
\[
\hyp43{1-n,1-n,2,2}{2+n,2+n,1}{1}=
\dfrac{\Gamma(n+1)^2\Gamma(n+2)^2\Gamma(4n+1)}{2(4n-1)\Gamma(2n+1)^4}
\]
From \eqref{8.1a} we get after some simplification
\[
\dfrac{\Gamma(2n+1)^2}{\Gamma(n+2)^2\Gamma(n+1)^2}
\hyp43{2,2,1-n,1-n}{1,2+n,2+n}{1}=\dfrac{\Gamma(4n-1)}{\Gamma(2n)\Gamma(2n+1)}
\]
The left-hand side is the starting summation. The right-hand side is the Catalan number $C_{2n-1}$. This proves the theorem. $\square$

\

\textbf{\fontsize{10.5}{12.5}\selectfont Proof of identity 2.9.}

We repeat \eqref{2.9}
\begin{equation}
S_8(n)=\sum_{k=0}^n\binom{n+k}{k}\dfrac{1}{2^k}=2^n
\label{9.1}
\end{equation}
The summation can be written as a hypergeometric function.
\begin{align*}
S_8(n)&=\sum_{k=0}^n\binom{n+k}{k}\dfrac{1}{2^k}=
\sum_{k=0}^n\dfrac{\Gamma(n+1+k)}{\Gamma(k+1)\Gamma(n+1)}\dfrac{1}{2^k}=
\sum_{k=0}^n(n+1)_k\dfrac{1}{k!}\dfrac{1}{2^k} \\
&=\sum_{k=0}^\infty(n+1)_k\dfrac{1}{k!}\dfrac{1}{2^k}
-\sum_{k=n+1}^\infty(n+1)_k\dfrac{1}{k!}\dfrac{1}{2^k} \\
&=\hyp21{-,n+1}{-}{\dfrac{1}{2}}
-\sum_{k=0}^\infty(n+1)_{k+n+1}\dfrac{1}{(n+1+k)!}\dfrac{1}{2^{k+n+1}} \\
&=2^{n+1}-\dfrac{1}{2^{n+1}}\dfrac{\Gamma(2n+2)}{\Gamma(n+1)\Gamma(n+2)}
\sum_{k=0}^\infty\dfrac{(2n+2)_k(1)_k}{(n+2)_k}\dfrac{1}{k!}\left(\dfrac{1}{2}\right)^k \\
&=2^{n+1}-\hyp21{2n+2,1}{n+2}{\dfrac{1}{2}}
\end{align*} 
\cite[7.3.7.(5)]{4} gives
\[
\hyp21{a,b}{\dfrac{a+b+1}{2}}{\dfrac{1}{2}}=
\sqrt{\pi}\dfrac{\Gamma((a+b+1)/2)}{\Gamma((a+1)/2)\Gamma((b+1)/2)}
\]
Application gives
\[
S_8(n)=2^{n+1}-\dfrac{1}{2^{n+1}}\dfrac{\Gamma(2n+2)}{\Gamma(n+1)\Gamma\left(n+\dfrac{3}{2}\right)}\sqrt{\pi}
\]
Using the duplication formula for the Gamma function \eqref{1.3} proves the theorem. $\square$

\

\textbf{\fontsize{10.5}{12.5}\selectfont Proof of identity 2.10.}

We repeat \eqref{2.10}
\begin{equation}
S_9=\sum_{k=0}^\infty\dfrac{(-1)^k}{(2k+1)^2\binom{2k}{k}}=
\dfrac{\pi^2}{6}-3\ln^2\left(\dfrac{\sqrt{5}-1}{2}\right)
\label{10.1}
\end{equation}
The summation can be written as a hypergeometric function.
\begin{align}
S_9&=\sum_{k=0}^\infty\dfrac{(-1)^k}{(2k+1)^2\binom{2k}{k}}=
\dfrac{\sqrt{\pi}}{4}\sum_{k=0}^\infty\dfrac{\Gamma\left(k+\dfrac{1}{2}\right)\Gamma(k+1)}{\Gamma\left(k+\dfrac{3}{2}\right)^2}\left(-\dfrac{1}{4}\right)^k \nonumber \\
&=\sum_{k=0}^\infty\dfrac{\left(\dfrac{1}{2}\right)_k(1)_k(1)_k}{\left(\dfrac{3}{2}\right)_k\left(\dfrac{3}{2}\right)_k}\dfrac{1}{k!}\left(-\dfrac{1}{4}\right)^k=
\hyp32{\dfrac{1}{2},1,1}{\dfrac{3}{2},\dfrac{3}{2}}{-\dfrac{1}{4}}
\label{10.1a}
\end{align}
\cite[7.4.3.(13)]{4} gives
\[
\hyp32{\dfrac{1}{2},1,1}{\dfrac{3}{2},\dfrac{3}{2}}{-z}=
\dfrac{1-x^2}{2x}\big[Li_2(x)-Li_2(-x)\big] \qquad \text{with} \qquad z(1-x^2)^2=4x^2
\]
$Li_2(x)$ is the dilogarithm function \cite{15}. With $z=\dfrac{1}{4}$ we get
\begin{equation}
\hyp32{\dfrac{1}{2},1,1}{\dfrac{3}{2},\dfrac{3}{2}}{-\dfrac{1}{4}}=
2\left[(Li_2\left(\sqrt{5}-2\right)-Li_2\left(2-\sqrt{5}\right)\right]
\label{10.3}
\end{equation}
For the right-hand side we use \cite[1.27]{15}
\[
Li_2(x)-Li_2(y)=Li_2\left(\dfrac{y(1-x)}{x(1-y)}\right)-Li_2\left(\dfrac{y}{x}\right)
-Li_2\left(\dfrac{1-x}{1-y}\right)+\dfrac{\pi^2}{6}-\ln(x)\ln\left(\dfrac{1-x}{1-y}\right)
\]
Using $x=\sqrt{5}-2=\left(\dfrac{\sqrt{5}-1}{2}\right)^3$ and $y=2-\sqrt{5}$ gives
\[
Li_2(\sqrt{5}-2)-Li_2(2-\sqrt{5})=Li_2\left(\dfrac{1-\sqrt{5}}{2}\right)-Li_2\left(\dfrac{\sqrt{5}-1}{2}\right)-Li_2(-1)+\dfrac{\pi^2}{6}-3\ln^2\left(\dfrac{\sqrt{5}-1}{2}\right)
\]
The polylogarithms in the right-hand side are known.

\cite[1.21]{15} \qquad gives $\qquad Li_2\left(\dfrac{1-\sqrt{5}}{2}\right)=-\dfrac{\pi^2}{15}+\dfrac{1}{2}\ln^2\left(\dfrac{\sqrt{5}-1}{2}\right)$

\cite[1.20]{15} \qquad gives $\qquad
Li_2\left(\dfrac{\sqrt{5}-1}{2}\right)=\dfrac{\pi^2}{10}-\ln^2\left(\dfrac{\sqrt{5}-1}{2}\right)$

\cite[1.9]{15} \qquad\, gives $\qquad Li_2(-1)=-\dfrac{\pi^2}{12}$

\

Application gives
\begin{equation}
Li_2(\sqrt{5}-2)-Li_2(2-\sqrt{5})=\dfrac{\pi^2}{12}-\dfrac{3}{2}\ln^2\left(\dfrac{\sqrt{5}-1}{2}\right)
\label{10.4}
\end{equation}
Combination of \eqref{10.1a}, \eqref{10.3} and \eqref{10.4} proves the theorem. $\square$

\

From \eqref{10.1} and \eqref{10.1a} there follows
\[
\hyp32{\dfrac{1}{2},1,1}{\dfrac{3}{2},\dfrac{3}{2}}{-\dfrac{1}{4}}=
\dfrac{\pi^2}{6}-3\ln^2\left(\dfrac{\sqrt{5}-1}{2}\right)
\]
which can be added to any table of hypergeometric series.

\

\section{Conclusion}
In this paper we prove a number of combinatorial identities. In some cases we use an identity of a hypergeometric series from Prudnikov et al. \cite{4}. However we do not know how these identities are derived. If these identities are proved with data of Gould \cite{13} then our proofs of these identities are not proofs. The book of Prudnikov et al. is much younger then the book of Gould. So there is a possibility  that the data of Gould is used. 

The same is for the identities from Gould \cite{13}. It is not known how these identities are proved. He hints on general series transformations. If he uses hypergeometric functions with transformations, he does the same as the method in this paper.

However, in this paper I only show a connection between some combinatorial identities and the hypergeometric functions. To what extent this is evidence, I leave up to the reader.

\

One can also wonder whether proving combinatorial identities still makes sense because nowadays this can be done using well-known computer programs. However
writing this paper gives great fun.

\end{document}